\documentclass{article}
\usepackage{amsfonts}
\usepackage{amsthm}
\usepackage{latexsym}
\usepackage{amssymb}
\newtheorem{thm}{Theorem}

\newtheorem{cor}{Corollary}[thm]

\newtheorem{quest}{Question}
\newtheorem{defn}{Definition}
\newtheorem{example}{Example}
\usepackage{graphicx}
\usepackage{epic}

\def\stab{{\rm Stab}}
\def\invar{{\rm Invar}}
\def\aut_G{{\rm Sym}}
\def\aut{{\rm Aut}}
\def\id{{\it id}}
\def\det{{\rm Det}}

\begin{document}

\begin{center}
\LARGE {\bf \textsc{Determining sets, resolving sets, and the exchange property}} 
\end{center}

\begin{center} \today \end{center}

\begin{center}
\textsc{Debra L. Boutin} 

\textsc{Department of Mathematics}

\textsc{Hamilton College, Clinton, NY 13323}

\textsc {\sl dboutin@hamilton.edu}

\end{center}

\begin{abstract} A subset $U$ of vertices of a graph $G$ is called a {\it determining set} if every automorphism of $G$ is uniquely determined by its action on the vertices of $U$.  A subset $W$ is called a {\it resolving set} if every vertex in $G$ is uniquely determined by its distances to the vertices of $W$. Determining (resolving) sets are said to have the {\it exchange property} in $G$ if whenever $S$ and $R$ are minimal determining (resolving) sets for $G$ and $r\in R$, then there exists $s\in S$ so that $S-\{s\}\cup \{r\}$ is a minimal determining (resolving) set.   This work examines graph families in which these sets do, or do not, have the exchange property.   This paper shows that neither determining sets nor resolving sets have the exchange property in all graphs, but that both have the exchange property in trees.   It also gives an infinite graph family ($n$-wheels where $n\geq 8$) in which determining sets have the exchange property but resolving sets do not.  Further, this paper provides necessary and sufficient conditions for determining sets to have the exchange property in an outerplanar graph.\end{abstract}

\section{Introduction}

A set of vertices $S$ of a graph $G$ is called a {\it determining set} if every automorphism of $G$ is uniquely determined by its action on the vertices of $S$.   The minimum size of  a determining set is a measure of graph symmetry and the sets themselves are useful in studying problems involving graph automorphisms \cite{AB1, AB2}.   The determining set was introduced in \cite{B1} and independently introduced as a {\it fixing set} in \cite{EH}.   Some transitive graphs have a small determining set compared to their vertex set and their automorphism group; others do not. For example, the hypercube  $Q_n$ has $2^n$ vertices and $2^n\cdot n!$ automorphisms, but has a determining set of cardinality only $\lceil \log_2n\rceil +1$ \cite{B3}.  In contrast, $K_n$ has $n$ vertices and $n!$ automorphisms, but requires fully $n-1$ vertices in a determining set. \medskip

Notice that by definition, the images of the vertices in a determining set under the trivial automorphism uniquely determine the positions of the remaining vertices.  Thus a determining set not only uniquely identifies each automorphism, but also uniquely identifies each vertex in the graph by its graph properties and its relationship to the determining set.  Historically, a number of different sets (some with multiple names) have been defined to identify the vertices of a graph relative to the given set.   Let $S\subseteq V(G)$.  $S$ has been called  a  {\it resolving set} \cite{CEJO} or {\it locating set} \cite{S} if every vertex of $G$ is uniquely identified by its distances from the vertices of $S$.   It has been called a   {\it distinguishing set} \cite{Ba} if each vertex can be uniquely identified by its set of neighbors that are in $S$.  It has been called a  {\it locating dominating set} \cite{RS} or {\it beacon set} \cite{CSS} if none of the sets of neighbors that are in $S$ is empty (that is, if $S$ is also a dominating set).  Each of these sets is a determining set, but not conversely.\medskip 

Since a determining set is the most general of these sets, it can be smaller (sometimes much smaller) than  the others. In this paper we will focus on comparing determining sets and resolving sets.   For instance, a smallest determining set for $Q_8$ has size $4$ \cite{B3}, while a smallest resolving set has size $6$ \cite{CHMPPSW}.   The $n$-wheel $W_n$ for $n\geq 4$ has a smallest determining set  of size $2$ (since its automorphisms are the same as those for the subgraph $C_n$) while for $n\geq 7$ a smallest resolving set has size $\lfloor {2n+2\over 5} \rfloor$ \cite{BCPZ} (because most pairs of vertices are at distance $2$).   Further, the size of a determining set can be easier to compute. The size of a smallest resolving set  for $Q_n$ is  known only for small $n$ and it takes \lq\lq laborious computations" to find them \cite{CHMPPSW},  while  for all $n$, ${\rm Det}(Q_n) = \lceil \log_2 n\rceil +1$\  \cite{B3}.\medskip

Both determining sets and resolving sets behave like bases  in a vector space in that each vertex in the graph can be uniquely identified relative to the vertices of these sets.  In fact, if a resolving set has minimum size then it is frequently called a {\it metric basis}  \cite{CHMPPSW} or just a {\it basis} \cite{CEJO} for the graph.  But though resolving sets and determining sets  do share some of the properties of bases, we will see in this paper that they do not necessarily share all the properties.  In particular, they do not always have the {\it exchange property}.  The form that the exchange property takes here is easily recognizable from linear algebra:  Determining (resolving) sets have the exchange property in $G$ if whenever $S$ and $R$ are minimal determining (resolving) sets and $r\in R$ then there exists $s\in S$ so that $S-\{s\}\cup\{r\}$ is a minimal determining (resolving) set.  In Section \ref{sect:basics} we will briefly discuss the connection with matroid theory.\medskip

There is a significant advantage if determining sets (or resolving sets) have the exchange property in a given graph: if the exchange property holds, then every minimal determining (resolving) set for that graph has the same size.  This makes algorithmic methods for finding the minimum size of such a set more feasible.\medskip
 
This paper examines graphs and graph families in which determining or resolving sets do (or do not) have the exchange property and is organized as follows.  Precise definitions for determining sets, resolving sets, and the exchange property are given in Section \ref{sect:basics}.  In Section \ref{sect:nothold} examples are given to show that neither determining nor resolving sets have the exchange property in all graphs.  Also given in this section are examples to show that, when the exchange property does not hold, a graph may have minimal determining (resolving) sets of different size.  Section \ref{sect:trees} gives criteria for a set to be a minimal determining (resolving) set in a tree and uses those criteria to show that determining (resolving) sets have the exchange property in trees.    Section \ref{sect:wheel} shows that there is an infinite family of graphs in which determining sets have the exchange property but resolving sets do not.  This infinite family is the set of all  $n$-wheels with $n\geq 8$.  Section  \ref{sect:outerplanar}  shows that neither resolving sets nor determining sets  have the exchange property in all outerplanar graphs.  In addition, it gives  necessary and sufficient conditions for determining sets to have the exchange property in in an outerplanar graph. Open questions are given in Section \ref{sect:quest}.

\section{Basics}\label{sect:basics}

Determining sets and resolving sets have already been introduced; the formal definitions follow. 

\begin{defn}\rm A subset $S$ of the vertices of a graph $G$ is called a {\it determining set} if whenever $g, h \in \aut(G)$  so that $g(s)=h(s)$ for all $s\in S$, then $g=h$.  The {\it determining number} of $G$, denoted  $\det(G)$, is the smallest integer $r$ so that $G$ has a determining set of size $r$.\end{defn}

The determining set has also been called the {\it fixing set} \cite{EH} and is an example of a {\it base of a permutation group action} \cite{DM}.\medskip

We need to be able to tell when a given set is a determining set.  Frequently this is done by looking at stabilizers.  Recall that for any subset $S\subseteq V(G)$,  the pointwise stabilizer of $S$ is $\stab(S) = \{g\in \aut(G) \ | \ g(v)=v, \ \forall v\in S\}  = \cap_{v\in S} {\rm Stab}(v)$.  

\begin{thm}\label{thm:stabid} \rm \cite{B1} A subset $S$ of the vertices of a graph $G$ is a determining set  if and only if $ \stab(S) = \{\id\}$. \end{thm}

Denote the standard distance between vertices $u$ and $v$ in a connected graph by $d(u,v)$.

\begin{defn}\rm A subset $S$ of the vertices of a connected graph $G$ is called a {\it resolving set} if for every $u,v\in V(G)$ there is $s\in S$ so that $d(s,v) \ne d(s,u)$. \end{defn}

We wish to study the exchange property with respect to determining sets and resolving sets. The exchange property is usually seen in connection with bases of a vector space, or more generally in matroid theory.  The following shows explicitly the connection with matroids.  \medskip

Define {\it d-independence} so that a set $S$ of vertices in a graph $G$ is d-independent if for every $s\in S$, $S-\{s\}$ is not a determining set.  With this definition, a maximal d-independent set is a minimal determining set.  One can similarly define  {\it r-independence} so that a maximal r-independent set is a minimal resolving set.  Each of these definitions of independence defines a hereditary system in the graph $G$.  Thus the question of whether the exchange property holds in $G$ is equivalent to the question of whether the hereditary system  in $G$ is a matroid \cite{W}.  Recall that one definition of the exchange property for hereditary systems is that whenever $S$ and $R$ are maximal independent sets and $r\in R$ then there exists $s\in S$ so that $S-\{s\}\cup\{r\}$ is a maximal independent set. (This form is  dual to the base exchange property given in \cite{W}.)  However, rather than saying that the \lq\lq hereditary system under d-independence (r-independence) has the exchange property in $G$" we will lazily say that determining sets (resolving sets) have the exchange property in $G$. More formally:

\begin{defn} \rm Determining (resolving) sets are said to have the {\it exchange property} in graph $G$ if whenever $S$ and $R$ are minimal determining (resolving) sets for $G$ and $r\in R$ then there is $s\in S$ so that $S-\{s\}\cup\{r\}$ is a minimal determining (resolving) set. \end{defn}

Recall that if the exchange property holds in a hereditary system then all maximal independent sets have the same size. Thus to show that the exchange property does not hold in a given graph, it is sufficient to show two minimal determining (resolving) sets of different size.   However, since the converse is not true, knowing that the exchange property does not hold does not guarantee that there are minimal determining (resolving) sets of different size.   

\section{Exchange Property Does Not Always Hold}\label{sect:nothold}

Recall that for $k<{n\over 2}$, the Kneser graph $K_{n:k}$ has vertices associated with the $k$-subsets of $\{1,\ldots, n\}$, and has vertices adjacent when their associated subsets are disjoint.  Here we will show that the exchange property fails for determining sets in $K_{7:3}$ and for resolving sets in $K_{5:2}$ (the Petersen graph).   However, we will also see that the exchange property holds for determining sets in $K_{5:2}$.

\begin{example} \rm In $K_{5:2}$, let $S=\{s_1=\{1,2\}, s_2=\{1,3\}, s_3=\{1,5\}\}$.  It is straightforward to show that $S$ is a minimal resolving set.  Let $r=\{2,4\}$.  Since $K_{5:2}$ is transitive, $r$ is an element of a minimal resolving set.  However for each $i=1,2,3$ it is straightforward to show that $S-\{s_i\}\cup  \{r\}$ is not a resolving set.  \end{example}

Further $K_{5:2}$ has minimal resolving sets of different size.  Both $S_1=$ $\{\{1,2\},$ $\{1,3\}, \{2,4\}, \{3,5\}\}$ and $S_2=\{\{1,2\}, \{1,3\}, \{1,5\}\}$ can be shown to be minimal resolving sets. \medskip

Thus the exchange property does not always hold for resolving sets, and further, minimal resolving sets for a given graph do not always have the same size.

\begin{example} \rm In $K_{7:3}$, let  $S=\{s_1=\{1,2,3\}, s_2=\{3,4,5\}, s_3=\{1,5,6\}\}$.  It is easy to show that $\stab(S)=\{\id\}$ and thus that $S$ is a determining set. $S$ can also be shown to be minimal.  Let $r=\{1,2,6\}$.  Since $K_{7:3}$ is transitive, $r$ is an element of a minimal determining set.  However for each $i=1,2,3$ it is straightforward to show that $S-\{s_i\}\cup \{r\}$ has nontrivial stabilizer and thus is not a determining set.  \end{example}

Further $K_{7:3}$ has minimal determining sets of different size.  Both $S_1=\{\{1,2,3\}, \{3,4,5\}, \{1,2,6\}, \{1,3,5\}\}$ and $S_2=\{\{1,2,3\}, \{3,4,5\}, \{1,5,6\}\}$ can be shown to be minimal determining sets. \medskip

Thus the exchange property does not always hold for determining sets and further, minimal determining sets for a given graph do not always have the same size.\medskip

However, the exchange property does hold for determining sets in $K_{5:2}$.  By brute force one can show that there are three different isomorphism classes of minimal determining sets. These all have size three.   For each isomorphism class one can show that an arbitrary added vertex can be exchanged for one of the original to yield a new minimal determining set. This gives the result.

\section{The Exchange Property in Trees}\label{sect:trees}

In this section we will see that both determining sets and resolving sets have the exchange property in trees.  In each case we start with a theorem giving a complete characterization of a determining or resolving set in a tree and then use the characterization to show that the exchange property holds.

\begin{defn} \rm Let $T$ be a tree.  Let $v\in V(T)$.  A {\it branch} at $v$ is a connected component of $T-\{v\}$ unioned with $v$.  Note that a branch is a maximal subgraph of $T$ containing $v$ as a leaf.   If $v$ is a vertex of degree at least $3$, a {\it branch path} at $v$ is a branch of the tree at $v$ that is a path.   \end{defn}

\begin{defn} \rm  Call vertex $v$ of $T$ an {\it exterior major vertex} if it has degree at least three and has at least one branch path.  \end{defn}

Slater provides  a nice characterization of a minimal resolving set of a tree as follows.

\begin{thm}\rm \cite{S}\label{treecriteria} A set $S$ of vertices is a minimal resolving set for a tree if and only if  for each major exterior vertex $v$ there exists a single vertex (different from $v$) in $S$ from precisely all but one of the branch paths at $v$.\end{thm}

Using this characterization we can prove the following theorem.

\begin{thm} \rm The exchange property holds for resolving sets in trees. \end{thm}

\begin{proof} Let $S$ be a minimal resolving set in a tree $T$.  Then for each major exterior vertex of $T$, $S$ contains a vertex from all but one of its branch paths.   Let $r$ be  an element of a minimal resolving set $R$.  Then $r$ is a vertex on a branch path, say $H$, of an exterior major vertex $v$ that has at least two branch paths.   If $S$ contains a vertex $s$ from $H$, then $S-\{s\}\cup\{r\}$  meets the requirements of Theorem \ref{treecriteria} and therefore is a minimal resolving set.  If $S$ does not contain a vertex from $H$, then $H$ is the one branch path at $v$ that does not provide such a vertex.  Since there are at least two block paths at $v$, there is some block path $K$ at $v$ (different from $H$) so that $S$ contains a vertex $s$ from $K$. Then $S-\{s\}\cup\{r\}$ meets the requirements of Theorem \ref{treecriteria} and  is therefore a minimal resolving set.\end{proof}

Finding a nice characterization of a minimal determining set in a tree will make it easy to prove that determining sets also have the exchange property in trees.  We begin the search for a characterization below.\medskip

Let $T$ be a tree with two adjacent vertices $x$ and $y$ in its center.  Let $T'$ be the result of adding a vertex $z$ of degree two to the edge between $x$ and $y$.  Then $z$ is the center of $T'$ and as such is fixed by all automorphisms of $T'$.  Further, $T$ and $T'$ have the same automorphism group (as permutations of $V(T) = V(T')-\{z\}$), and all vertex stabilizers are the same.  The vertex $z$ is redundant in any determining set of $T'$ since $\stab(z) = \aut(T') = \aut(T)$.   Thus every determining set for $T$ is a determining set for $T'$ and every determining set for $T'$ that does not include $z$ is a determining set for $T$.  In particular, $T$ and $T'$ have the same minimal determining sets. \medskip

Thus for the purposes of studying minimal determining sets, we may assume that a tree has a single vertex as its center.  

\begin{thm}\label{char}\rm Let $T$ be a tree with a single vertex as its center.  A set $S$ of vertices is a determining set for $T$ if and only if for each vertex $v$ there exists a vertex (different from $v$) in $S$ from all but one of the branches from each isomorphism class of branches at $v$. \end{thm}

\begin{proof}

\noindent $\Longrightarrow$  \ \  Suppose there is a vertex $v$ with two isomorphic branches $H_1$ and $H_2$ so that no vertex (different from $v$) of either branch is in $S$.  Then there is an automorphism $\alpha$  that transposes $H_1$ and $H_2$,  fixes $v$, and fixes all vertices outside of $H_1$ and $H_2$.  Then  $\alpha$ is nontrivial and $\alpha\in \stab(S)$.  Thus $S$ is not a determining set.\medskip

\noindent $\Longleftarrow$  \ \  Suppose that $S$ is not a determining set.  Then there is some nontrivial $\alpha\in \stab(S)$. Thus there exists $x\in V(T)$ so that $\alpha(x)\ne x$.  Let $v$ be the first vertex on a path from $x$ to the center of $T$ that is fixed by $\alpha$.  (Since $T$ has a single vertex at its center, such a $v$ must exist.) Let $H$ be the  branch of $v$ containing $x$. Since $\alpha$ moves every vertex of $H$ except $v$, no vertex of $H$ different from $v$  can be in $S$.   Note that $\alpha(H)$ is also a branch at $v$ and is isomorphic to and distinct (except for $v$) from $H$. Again, since $\alpha$ moves every vertex of $\alpha(H)$ except $v$, no vertex of $\alpha(H)$ different from $v$  can be in $S$.  Thus $H$ and $\alpha(H)$ are isomorphic branches at $v$ and neither contributes a vertex to the set $S$.\end{proof}

The example below shows that the hypothesis that $T$ has a unique center is necessary in the theorem above. 

\begin{example} \rm Let $T$ be the graph in Figure \ref{tree}. Vertices $x$ and $y$ are the center of $T$.  Note that no vertex has a nontrivial branch isomorphism class.  Thus if the theorem held for trees with non-unique centers, this graph would have the empty set as a determining set.  However, there is an automorphism that transposes the center vertices and their branches that is not captured by this (empty) set.  As discussed earlier, we could place a vertex between $x$ and $y$, and this would be the new center. We see that this new vertex would have a nontrivial branch isomorphism class and that there would be a vertex from one of these branches in any determining set.  This would provide the correct result.\end{example}

\begin{figure}[htb]
\begin{center}
\thicklines
\setlength{\unitlength}{2.25mm}
\begin{picture}(35,7)

\put(0,2){\circle*{1}}

\put(5,2){\circle*{1}}

\put(10,2){\circle*{1}}

\put(15,2){\circle*{1}}

\put(20,2){\circle*{1}}

\put(25,2){\circle*{1}}

\put(30,2){\circle*{1}}

\put(35,2){\circle*{1}}

\put(10,7){\circle*{1}}

\put(25,7){\circle*{1}}

\put(14,-1){$x$}

\put(19,-1){$y$}

\put(0,0){\drawline(0,2)(5,2)}
\put(0,0){\drawline(5,2)(10,2)}
\put(0,0){\drawline(10,2)(15,2)}
 
\put(0,0){\drawline(15,2)(20,2)}
\put(0,0){\drawline(20,2)(25,2)}
\put(0,0){\drawline(25,2)(30,2)}
\put(0,0){\drawline(30,2)(35,2)}

\put(0,0){\drawline(25,2)(25,7)}
\put(0,0){\drawline(10,7)(10,2)}
\end{picture}
\end{center}
\caption{} \label{tree}
\end{figure}

\begin{thm}\rm The exchange property holds for determining sets in trees.\end{thm}

\begin{proof} If $\det(T)=0$ then the result holds vacuously.  So assume that $\det(T)>0$.\medskip

Recall that if the tree $T$ does not have a unique center, we can add a vertex of degree $2$ to the edge between the two vertices of the center, yielding a new tree $T'$ with a unique center.  Since $T'$ has precisely the same minimal determining sets as $T$, determining sets have the exchange property in $T'$ if and only if they have the exchange property in $T$. Thus, we may assume $T$ has a unique center.\medskip

Let $S$ be a minimal determining set for $T$.  Let $R$ be another and $r\in R$. Let $v$ be the first vertex on the path from $r$ to the center of $T$ with the property that if $H$ is the branch at $v$ containing $r$, then there exists a branch $K$ at $v$ isomorphic to $H$.  That is, $r$ is a vertex on one of the branches in a nontrivial isomorphism class of branches at $v$.  Note that since $R$ is a nonempty minimal determining set, such a $v$ must exist.\medskip

If there is $s\in S$ so that $s\ne v$ and $s$ is on $H$, then exchanging $s$ for $r$ yields a determining set because $S'=S-\{s\}\cup \{r\}$ continues to meet the conditions of Theorem  \ref{char}. If there is no such $s$ on $H$, then $H$ is the one branch in its isomorphism class of branches at $v$ that does not contain an element of $S$.  Then there is some $s\in S$ so that $s\ne v$ and is on $K$.   Again, exchanging $s$ for $r$ yields a determining set.\medskip

It is intuitively clear that the minimality of $S$ implies the minimality of $S'$.  However, set theoretic arguments using vertex stabilizers can be used to prove this more formally.   The proofs would be similar to (but slightly simpler than) the proof for minimality in Case 2 in Theorem \ref{OPEP}.
 \end{proof}


\section{Wheels}\label{sect:wheel}

So far the exchange property results for determining sets and for resolving sets do not look very different.  Neither set always has the exchange property, and both have the exchange property in trees.  However, in this section we see that the set of $n$-wheels where $n\geq 8$ is an infinite family of graphs in which the exchange property holds for determining sets but not for resolving sets.\medskip

Recall the the $n$-wheel  is the join of $C_n$ and $K_1$ for $n\geq 3$.  That is, $W_n$ contains an $n$-cycle with a single additional vertex adjacent to all the vertices of the cycle.

\begin{thm} \rm Determining sets have the exchange property in $n$-wheels. \end{thm}

\begin{proof} Recall that $W_3=K_4$.  It is trivial to see that determining sets have the exchange property in complete graphs.  When $n\geq 4$,  $C_n$ and $W_n$ have the same automorphism group and the same vertex stabilizers on the vertices of the $n$-cycle.  Thus the minimal determining sets for $W_n$ are precisely the minimal determining sets for $C_n$, and these are the pairs of non-antipodal vertices on the $n$-cycle. If we have a set $S=\{u,v\}$ in $V(W_n)$ of non-antipodal vertices and another vertex $x\ne u,v$ of a minimal determining set ({\it i.e.}, $x$ is not the central vertex of $W_n$), then $x$ is non-antipodal to at least one of $u$ and $v$ on the cycle and can be exchanged for the other. Thus determining sets have the exchange property in all $n$-wheels.\end{proof}

\begin{thm} \rm For $n\geq 8$, resolving sets do not have the exchange property in $W_n$. \end{thm}

\begin{proof} Label the vertices of degree three in $W_n$ by $1, 2, \ldots, n$ around the cycle.  Buczkoski, Chartrand, Poisson, and Zhang \cite{BCPZ} give criteria for a set of vertices in $W_n$ to be a resolving set. The criteria use the size of the {\it gaps}, the distance along the cycle between sequentially ordered vertices in the set.  A set $S$ is a resolving set if and only if there is no gap of size four or greater, there is at most one gap of size three, and any gap of size greater than one must have both neighboring gaps of size one or less. Further, Buczkoski {\it et al.} give minimal resolving set for each $n\geq 7$ depending on its residue modulo $5$.  These, or slight modifications of these, are used in the following. 

\begin{itemize}  

\item Suppose that $n=5k$ for $k\geq 2$.  Then $S=\{5i+1, 5i+4: 0\leq i \leq k-1\}$ is a minimal resolving set. There is no $s\in S$ so that $S-\{s\}\cup \{2\}$ is a resolving set.  In particular the only  possibilities are the elements of $S$ nearest to $2$.  However, removing $1$ would yield neighboring gaps $\{n,1\}$ and $\{2,3\}$, while removing $4$ would leave neighboring gaps $\{5,6\}$ and $\{7,8\}$.

\item Suppose that $n=5k+1$ for $k\geq 2$.  Then $S=\{5i+1, 5i+4: 0\leq i\leq k-2\}\cup\{n-3, n-1\}$ is a minimal resolving set.  There is no $s\in S$ so that $S-\{s\}\cup \{2\}$ is a resolving set.  In particular the only  possibilities are the elements of $S$ nearest to $2$.  However, removing $1$ would yield neighboring gaps $\{n, 1\}$ and $\{n-6,n-5,n-4\}$, while removing $4$ would leave gap $\{3, 4, 5\}$ a second gap of size three. 

\item Suppose that $n=5k+2$ for $k\geq 1$.  Then $S=\{1, 5\} \cup \{5i+2, 5(i+1)): 1\leq i\leq k-1\}\cup\{n\}$ is a minimal resolving set.  (This is different from the minimal resolving set given in \cite{BCPZ}.)   There is no $s\in S$ so that $S-\{s\}\cup \{6\}$ is a resolving set.  In particular, removing $5$ would yield a gap $\{2,3,4,5\}$ of size larger  than three, while removing $7$ would leave two gaps $\{2,3,4\}$ and $\{8,9,10\}$ of size three.

\item Suppose that $n=5k+3$ for $k\geq 2$.  Then $S=\{ 5i+1, 5i+4: 0\leq i\leq k-2\}\cup\{n-7, n-3, n-1\}$ is a minimal resolving set for $W_n$.  There is no $s\in S$ so that $S-\{s\}\cup\{7\}$ is a resolving set.  In particular, removing $6$ would yield neighboring gaps $\{5,6\}$ and $\{8,9\}$, while removing $10$ would leave gap $\{8,9,10,11\}$.  Suppose that $n=8$.  Then $S=\{ 1, 5, 7\}$ is a minimal resolving set.  Further there is no $s\in S$ so that $S-\{s\}\cup \{8\}$ is a resolving set.

\item Suppose that $n=5k+4$ for $k\geq 1$. Then $S=\{ 5i+1, 5i+4: 0\leq i\leq k\}$is a minimal resolving set for $W_n$.  There is no $s\in S$ so that $S-\{s\}\cup\{3\}$ is a resolving set. In particular, removing $1$ would yield neighboring gaps $\{n-2, n-1\}$ and $\{1,2\}$, while removing $4$ would leave neighboring gaps $\{4,5\}$ and $\{7,8\}$.\end{itemize} \end{proof} 


\section{Outerplanar Graphs}\label{sect:outerplanar}

We will see in the examples below that neither determining sets nor resolving sets have the exchange property in all outerplanar graphs.   However, it is straightforward to tell whether or not the exchange property holds for determining sets in a given outerplanar graph.  We first find criteria for the exchange property to hold for determining sets in a $2$-connected outerplanar graph. Along the way we find a description for determining sets in such a graph.  Finally, we use these to prove analogous results for determining sets in a connected outerplanar graph. 

\begin{example}\label{resolvingexample}\rm   Let $G$ be the $8$-wheel with outer vertices labeled $1, \ldots, 8$ and the edge from $1$ to $8$ removed.  This graph is outerplanar.  Each of the sets $\{1,5,7\}$ and $\{1, 3, 6, 8\}$ can be shown to be a minimal resolving set for $G$.  Thus resolving sets do not have the exchange property in $G$.\end{example}

\begin{example}\label{OP2Cexample} \rm Consider the graph in Figure \ref{properdihedral}.  Note that $\{x\}$ (or any single vertex that is not on the $3$-cycle) has a trivial stabilizer and thus forms a minimal determining set.  Also note that $u$ and $v$  (or any pair of vertices from the $3$-cycle) are fixed by distinct reflections of the Hamilton cycle.  Then $\stab(\{u,v\})$ is trivial and thus  $\{u,v\}$ is also a minimal determining set.  Since there are minimal determining sets of different size, the  exchange property does not hold for determining sets in this graph.

\begin{figure}[htb]
\centerline {
\scalebox{1.0} {\includegraphics[width=1.5in]{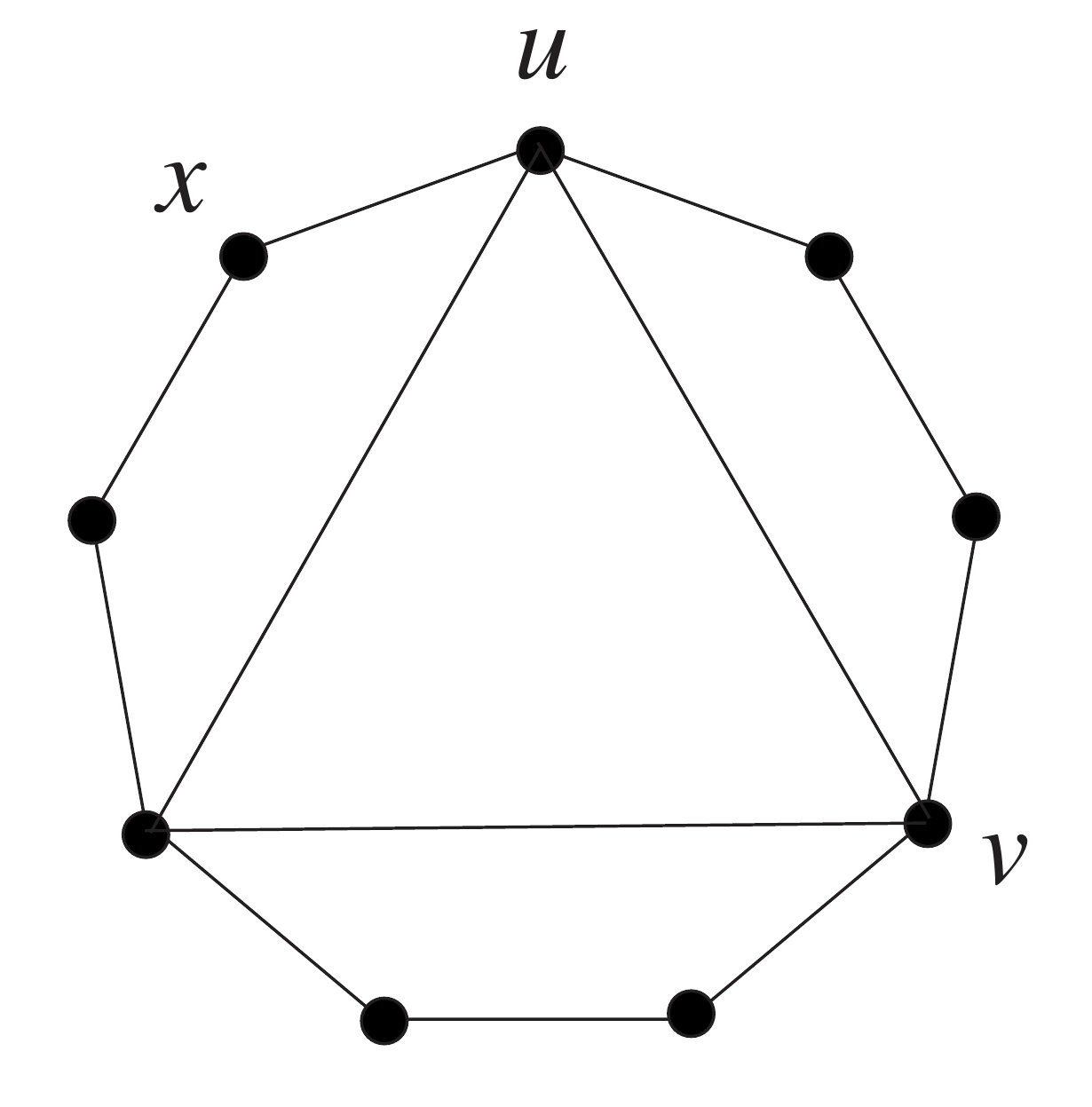}}
}
\caption{} \label{properdihedral}
\end{figure}\end{example}

Recall that a $2$-connected outerplanar graph with at least three vertices contains a unique Hamilton cycle \cite{CB}.  The uniqueness of the cycle guarantees that it is preserved by every graph automorphism.  Thus for such a graph $G$, $\aut(G)$ is a subgroup of $D_{|V(G)|}$, the automorphism group of the cycle.  That is, the automorphisms of $G$ are realized by rotations or reflections of the Hamilton cycle of $G$.  When $G$ is $K_2$, its automorphism group contains a single nontrivial reflection (or a single nontrivial rotation). Here we will only consider reflections that are nontrivial.

\begin{thm}\label{OP2CEP} \rm  Determining sets fail to have the exchange property in a $2$-connected outerplanar graph $G$ if and only if $G$ contains a vertex moved by every reflection and contains two vertices fixed by different reflections.\end{thm}

\begin{proof}  If $\aut(G)$ is trivial, the result is vacuously true.  Suppose that $G$ is a $2$-connected graph with nontrivial automorphisms.   Recall that no vertex of $G$ is fixed under any rotation of the Hamilton cycle.  Thus each vertex stabilizer in $G$ is either trivial or is generated by a single reflection.\medskip

\noindent $\Longrightarrow$ \ \  Suppose that $G$ has at least one vertex that is moved by every reflection, but no two vertices that are fixed by different reflections. (This includes the case in which there are no reflections in $\aut(G)$.) Any vertex that is moved by every  reflection has a trivial stabilizer and thus forms a (minimal) determining set.  Since no two vertices are fixed by different reflections, but each nontrivial stabilizer is generated by a single reflection, all nontrivial stabilizers are equal.  Then the intersection of multiple vertex stabilizers is trivial if and only if one of the stabilizers is trivial.  Thus a minimal determining set must be composed of a single vertex that is moved by every reflection.  Clearly, two such minimal determining sets can have their vertices exchanged.  Thus the exchange property holds in this case.\medskip

Suppose that $G$ has no vertex moved by every reflection and has two vertices that are fixed by different reflections.  The intersection of the vertex stabilizers of two vertices fixed by different reflections is trivial and thus these vertices form a (minimal) determining set.  Since every vertex is fixed by some reflection, no vertex stabilizer is trivial.  Thus the presence of two vertices fixed by different reflections is both necessary and sufficient for a determining set.  Suppose $S=\{s_1,s_2\}$ and $R=\{r_1,r_2\}$ are two such minimal determining sets.  Choose $r_1\in R$.  Since $s_1$ and $s_2$ are fixed by different reflections, one of them is fixed by a reflection other than the one that fixes $r_1$.  Suppose that $s_2$ is such a vertex.  Then $S-\{s_1\}\cup\{r_1\}=\{r_1,s_2\}$ is a minimal determining set for $G$. Thus the exchange property holds in this case.\medskip

\noindent $\Longleftarrow$ \ \ Suppose that $G$ has a vertex, say $x$, that is moved by every  reflection, as well as vertices, say $u$ and $v$, that are fixed by different reflections.  Then $\{x\}$ and $\{u,v\}$ are  minimal determining sets of different size.  Thus the exchange property does not hold for determining sets in $G$. Note that such a graph is displayed in Figure \ref{properdihedral}  (Example \ref{OP2Cexample}). \end{proof}

The proof of Theorem \ref{OP2CEP} yields the following corollary.

\begin{cor}\label{OP2Ccriteria}  \rm Let $G$ be a $2$-connected outerplanar graph with nontrivial automorphisms.  If $G$ contains a vertex $x$ that is moved by every  reflection, then $\{x\}$ is a minimal determining set.  If $G$ contains two vertices $u$ and $v$, each fixed by a different reflection, then $\{u,v\}$ is a minimal determining set. These are the only possible minimal determining sets for $G$.  \end{cor}
Before heading into the theorems regarding general outerplanar graphs, it will be useful to collect ideas about the structure of these graphs and facts regarding their automorphisms.\medskip

Each outerplanar graph $G$ has an associated nontrivial block-cutvertex tree \cite{W}.  Denote this by $T$.  That is, $T$ is a tree with a  $B$-vertex for each $2$-connected block of $G$ and a $C$-vertex for each cutvertex of $G$, and edges between $B$-vertices and their associated $C$-vertices.  Call a subgraph of $G$ a {\it blockbranch} if its associated block-cutvertex subtree is a branch of  $T$ attaching at a $C$-vertex.   There are times when it is useful to consider a {\it trivial blockbranch}; this is simply a vertex of the block.  Note that if $G$ is $2$-connected  the block-cutvertex tree consists of a single $B$-vertex.\medskip

Any path from a $C$-vertex to a leaf must be of odd length, and any path from a $B$-vertex to a leaf must be  of even length. Thus the center of the block-cutvertex tree cannot contain both a $C$-vertex and a $B$-vertex.  Further since no two $B$- or $C$-vertices are adjacent, the center cannot contain two $B$-vertices or two $C$-vertices.  Thus,  the center of a block-cutvertex tree is a unique vertex, either a $B$-vertex or a $C$-vertex. If the center of $T$ is a $B$-vertex, call the associated block of $G$ the {\it central block}. Otherwise a block is called a {\it non-central block}. If the center of $T$ is a $C$-vertex, call the associated cutvertex the {\it central cutvertex}.  It will ease our notation to refer to the vertices of the central block or the central cutvertex as $Z$.   Note that $Z$ does not necessarily contain the center of $G$.\medskip

Since graph automorphisms take $2$-connected blocks to $2$-connected blocks, and cutvertices to cutvertices, they preserve the block structure of $G$.  In particular,  if $v$ is a vertex of a block $B$, $\alpha\in \aut(G)$, and $\alpha(v)$ is also in $B$, then $B$ is invariant under $\alpha$.  If $v$ is a vertex (other than the terminal cutvertex) of a blockbranch $H$, and $\alpha(v)$ is also in $H$, then $H$ is invariant under $\alpha$.  Moreover, every $\alpha\in\aut(G)$ induces an automorphism, $\tilde \alpha$, on the associated block-cutvertex tree.   Since every automorphism  of a tree preserves its center, every automorphism of $G$ fixes $Z$ setwise.\medskip

If a block $B$ of $G$ is non-central then there is a cutvertex of $B$  that is closest to $Z$. This vertex is unique among the vertices of $B$ and thus must be fixed by any automorphism under which $B$ is invariant.  In particular, if a non-central block is invariant but not fixed under an automorphism, then the block has at least three vertices and the action of the automorphism on the Hamilton cycle is a reflection through the cutvertex of $B$ that is closest to $Z$, with analogous action on the attached blockbranches.\medskip

For any block $B$ of $G$, we can define a {\it block determining set} for $B$ to be a subset $U$ of $V(B)$ with the property that any automorphism of $G$ that fixes $U$ pointwise must also fix $B$ pointwise. More formally, let $\invar(B)$ be the set of automorphisms of $G$ under which $B$ is invariant.  It is easy to see that $\invar(B)$ is a group and that $\stab(B) \triangleleft \invar(B)$.  The symmetries of $B$ under the automorphisms of $\aut(G)$ are captured precisely by the quotient group $\invar(B) \slash \stab(B)$, which we will denote by $\aut_G(B)$.  Note that $\aut_G(B)$ is a  (often proper) subgroup of the automorphisms of $B$ when thought of as a graph in its own right.  We see that  a subset $U$ of $B$ is a block determining set for $B$  if and only if the (pointwise) stabilizer of $U$ in $\aut_G(B)$ is trivial.  Note that if $B$ is non-central and $\aut_G(B)$ is nontrivial then the nontrivial symmetry of $B$ in $\aut_G(B)$ is a reflection.  Therefore a minimal block determining set for the block is formed by any vertex of $B$ that is moved by the reflection in $\aut_G(B)$. If $B$ is the central block, its block determining sets follow the criteria of Corollary \ref{OP2Ccriteria} under  the group $\aut_G(B)$.\medskip

We now proceed to state criteria for a set of vertices of an outerplanar graph to be a determining set.

\begin{thm}\label{OPcriteria} \rm  Let $G$ be a connected outerplanar graph and let $T$ be its block-cutvertex tree. A set of vertices $S$ is a determining set for $G$ if and only if each of the following conditions hold.
\begin{enumerate}  

\item If $B$ is a block of $G$ then $X=\{v\in V(B)  \ | \ S$ contains a vertex of a blockbranch of $B$ at $v\}$ is a block determining set for $B$.  

\item If $v$ is a cutvertex of $G$ then $S$ contains a vertex (other than $v$) from all but one of the blockbranches in each isomorphism class of blockbranches at $v$. 
 
 \end{enumerate}  \end{thm}
 
\begin{proof} 

Let $S$ be a subset of vertices of $G$. Assume that $S$ is not a determining set for $G$.  Then there exist some nontrivial $\alpha\in \stab(S)$.\medskip

Suppose there is a block $B$ of $G$ that is invariant but not fixed under $\alpha$.  Since $\alpha\in \stab(S)$ leaves invariant all blockbranches of $B$ that contain elements of $S$, $\alpha$ fixes in $B$ the attaching cutvertices of these blockbranches. That is, $\alpha$ fixes the vertices of $X$.  But since $\alpha$ does not fix $B$, this means $X$ is not a block determining set for $B$.  This violates Condition 1.\medskip

Suppose no block is invariant but not fixed under $\alpha$.  Since $\alpha$ is nontrivial, this means that the induced automorphism $\tilde \alpha$ is nontrivial on $T$.  Then $\tilde \alpha$ permutes branches from a common vertex $\tilde y$ of $T$.  That is, there are isomorphic branches $\tilde H$ and $\tilde \alpha(\tilde H)$ at $\tilde y$, which are disjoint except for their common cutvertex. Thus in $G$, there are associated blockbranches $H$ and $\alpha(H)$ which either have a common cutvertex or attach to a common block at different cutvertices.  Suppose that $H$ and $\alpha(H)$ have different attaching cutvertices, say $v$ and $\alpha(v)$, on a common block $B$. Since $v, \alpha(v)$ are distinct and both in $B$, $B$ is invariant but not fixed under $\alpha$. This contradicts our assumption on $\alpha$.  Thus $H$ and $\alpha(H)$ share a common cutvertex $v$ and the only vertex of these blockbranches that is fixed by $\alpha$ is $v$.  But since $\alpha$ fixes all the vertices in $S$ this means that no vertex of $H$ or $\alpha(H)$ other than $v$ can be contained in $S$.  This violates Condition 2.\end{proof}

The above criteria for a determining set for an outerplanar graph allows us to prove the following.

\begin{thm}\label{OPEP}\rm  The exchange property fails for determining sets in a connected outerplanar graph if and only if it has a central block that has trivial (pointwise) stabilizer, and that contains a vertex moved by every  reflection of $B$ under $\aut_G(B)$, as well as two vertices fixed by different reflections of $B$ under $\aut_G(B)$. \end{thm}

\begin{proof} \noindent $\Longrightarrow$  \ \ Let $G$ be a connected outerplanar graph that does not meet the hypothesis above.  That is, $G$ does not have a central block $B$ with trivial stabilizer that contains a vertex moved by every  reflection of $B$, as well as two vertices fixed by different reflections of $B$.\medskip

If $\det(G)=0$ then the result holds vacuously.  Assume that $\det(G)>0$.\medskip

Let $T$ be the block-cutvertex tree for $G$. Let $S$ and $R$ be minimal determining sets for $G$. Let $r\in R$. Since $R$ is a minimal determining set there is at least one $\alpha\in \aut(G)$ so that $\alpha(r) \ne r$. Let $\tilde x$ be a $B$-vertex of $T$ whose associated block contains $r$.  For each $\alpha$ that moves $r$, let $\tilde y$ be the first vertex in a path from $\tilde x$ to the center of $T$ that is fixed by $\tilde \alpha$.  Over all such $\alpha$, choose one for which $\tilde y$ is closest to $\tilde x$.  (Note that if $r$ and $\alpha(r)$ are in the same block then $\tilde x = \tilde y$.)  Denote the (possibly trivial) branch from $\tilde y$ containing $\tilde x$ by $\tilde H$ and its associated blockbranch in $G$ by $H$.\medskip

By our choice of $\tilde y$, if it is a $C$-vertex then $H$ and $\alpha(H)$ are disjoint except for a common cutvertex, say $v$, and $\alpha$ permutes blockbranches  isomorphic to $H$ at $v$.  If $\tilde y$ is a $B$-vertex then $H$ and $\alpha(H)$ have distinct cutvertices on a common block $B$, and $\alpha$ acts as a symmetry on the block $B$ while carrying along isomorphic block branches.  Further, by our choice of $\tilde y$, $r$ does not satisfy Condition 1 for any block, nor Condition 2 for any cutvertex, that lies strictly between $r$ and the block or cutvertex associated with $\tilde y$.  Since $R$ is a minimal determining set this means that $r$ is necessary to fulfill Condition 1 for the block associated with $\tilde y$ or Condition 2 for the cutvertex associated with $\tilde y$. No other element of $R$ fulfills this role. \medskip 

\noindent{\bf Case 1:}  Suppose that $\tilde y$ is a $B$-vertex with associated block $B$.  Then $B$ is invariant but not fixed under $\alpha$.  Since $R$ and $S$ are determining sets for $G$, by Theorem \ref{OPcriteria} each provides a block determining set for $B$.  That is, each of $X=\{v\in V(B) \ | \  R$ contains a vertex of a blockbranch of $B$ at $v\}$ and $Y=\{w\in V(B) \ | \  S$ contains a vertex of a blockbranch of $B$ at $w\}$ are block determining sets for $B$.  This means that the stabilizers of $X$ and $Y$ in $\aut_G(B)$ are trivial. Let $v$ be the element of $X$ at which the blockbranch of $B$ containing $r$ is attached. \medskip

Recall that Corollary \ref{OP2Ccriteria} extends to block determining sets of $B$ under $\aut(B)$.  In particular, each block determining set must either contain a vertex moved by every reflection of $B$ under $\aut_G(B)$, or contain two vertices that are fixed by different reflections of $B$ under $\aut_G(B)$.  This allows us to break Case 1 in three subcases.  In each subcase we will identify a vertex $w$ of $Y$ that performs the same function in a block determining set for $B$ as $v$ does.  In particular we will find $w\in Y$ so that $Y-\{w\}\cup \{v\}$ and $X-\{v\}\cup\{w\}$ are both block determining sets for $B$.  We will then use $w$ to identify a vertex $s\in S$ to exchange for $r\in R$.  Finally we will prove that this exchange yields a minimal determining set.\medskip

\noindent{\bf Case 1.1:} Suppose that every block determining set for $B$ contains a vertex moved by every reflection in $\aut_G(B)$. (This includes the case where there is no reflection in $\aut_G(B)$.)\medskip

Then each of $X$ and $Y$ contains a vertex with trivial stabilizer in $\aut_G(B)$.  By the minimality of $R$, vertex $v$ plays this role for $X$.  Suppose $w$ plays this role for $Y$.  Then $Y-\{w\}\cup \{v\}$ and $X-\{v\}\cup \{w\}$ both have trivial stabilizer in $\aut_G(B)$.\medskip

\noindent{\bf Case 1.2:} Suppose that every block determining set for $B$ contains two vertices fixed by different reflections in  $\aut_G(B)$.\medskip

Then there exist  $v_1,v_2\in X$ and $w_1, w_2\in Y$ so that the vertices in each pair are fixed by different reflections in $\aut_G(B)$. By the minimality of $R$, one of $v_1$ and $v_2$ is $v$.  Assume $v_1=v$.  If one of $w_1, w_2$ is fixed by the same reflection as $v$ then let $w$ be that vertex. In this case, since both $X$ and $Y$ have trivial stabilizers in $\aut_G(B)$, both $Y-\{w\}\cup\{v\}$ and $X-\{v\}\cup \{w\}$ do also.  Suppose that neither of $w_1,w_2$ is fixed by the same reflection as $v$.  Since $w_1$ and $w_2$ are fixed by different reflections, at least one of them must be fixed by a reflection different from that of $v_2$.  Let this be $w=w_1$.  Then by our choice of $w$, and since $X-\{v\}\cup\{w\}$ contains both $v_2$ and $w$, it has  trivial stabilizer in $\aut_G(B)$. Further, since by hypothesis $v$ is fixed by a reflection different from that of $w_2$, and since $Y-\{w\}\cup\{v\}$ contains both $w_2$ and $v$, it also has trivial stabilizer in $\aut_G(B)$.\medskip

\noindent{\bf Case 1.3:}  Suppose that a block determining set for $B$ can contain either a vertex moved by every reflection in $\aut_G(B)$ or two vertices fixed by different reflections in $\aut_G(B)$.\medskip

 By hypothesis on $G$, the above assumption means that the pointwise stabilizer of $B$ in $\aut(G)$ is nontrivial. Then there are elements of $X$ and $Y$ that are there by virtue of automorphisms in the stabilizer of $B$. That is, there are vertices of $R$ and $S$ fulfilling Condition 1 for non-central blocks or Condition 2 for cutvertices of $G$. (Note that since there are two different reflections in $\aut_G(B)$, $B$ must be the central block.) Let $U$ be the subset of $X$  consisting of vertices coming from blockbranches that are invariant but not fixed under automorphisms that fix $B$. Notice that $U$ is also a subset of $Y$, that $v$ is not in $U$, and that $U$ is nonempty by hypothesis.  If the stabilizer of $U$ in $\aut_G(B)$ is trivial then $v$ is unnecessary for $X$ to be a block determining set for $B$.  Then $r$ is not necessary to fulfill Condition 1 for $B$. But we have argued that the choice of $B$ and the minimality of $R$ prevent this situation. Thus the stabilizer of $U$ in $\aut_G(B)$ is not trivial, which means that it is generated by a single reflection.  Since $X$ and $Y$ are block determining sets for $B$, they each contain a vertex of $B$ that is not fixed by the reflection that fixes $U$.  The vertex $v$ fills this role for $X$; let $w$ denote a vertex that fills the role for $Y$.  Then $Y-\{w\} \cup \{v\}$ contains $U$ and $v$ and therefore has trivial stabilizer.  Similarly $X-\{v\}\cup \{w\}$ contains $U$ and $w$ and therefore has trivial stabilizer.\medskip

Using our choice of $w$ from the appropriate case above, let $s$ be a vertex of $S$ on a blockbranch of $B$ attached at $w$.  Let $S'=S-\{s\}\cup \{r\}$. We will use the fact that $X-\{v\}\cup\{w\}$ has trivial stabilizer in $\aut_G(B)$ to show that $\stab(Y-\{w\}\cup\{v\}) \subseteq \stab(s)$. We will then use this fact to show that $S'$ is a determining set.\medskip

Suppose that $\beta\in \stab(Y -\{w\} \cup \{v\})$.  Then since $Y'=Y-\{w\}\cup \{v\}$ is a block determining set for $B$, $\beta$ fixes $B$ and therefore $w$.  Thus the blockbranches of $B$ at $w$ are invariant under $\beta$. Suppose that $\beta(s)\ne s$. Then either $s$ fulfills Condition 1 of Theorem \ref{OPcriteria} for some block between $s$ and $w$, or $s$ fulfills Condition 2 for some cutvertex between $s$ and $w$ (including the possibility of $w$).  But since $R$ is a determining set, $R$ must also contain a vertex, say $t$,  of a blockbranch at $w$ fulfilling the same condition for the same block or cutvertex.   Thus $w$ is in $X$. Since $X-\{v\} \cup \{w\}$ is a block determining set for $B$, this means $v$ (and therefore $r$) is not necessary to fulfill Condition 1 for $B$.  But as we've argued, this situation cannot occur. Thus $\stab(Y')\subseteq\stab(s)$.\medskip

Let $W$ be the subset of vertices from $S$ on blockbranches of $B$ at elements of $Y$.  Let  $W'=W-\{s\}\cup \{r\}$.  Note that if $\alpha$ fixes the elements of $W$ (resp. $W'$) then it fixes the elements of $Y$ (resp. $Y'$).  Combining the with the results above yields  $\stab(W')\subseteq\stab(Y')\subseteq\stab(s)$. \medskip

Since $\stab(W') \subseteq \stab(s)$ and since $W-\{s\}\subset S$, basic set theory gives us that $\stab(S') = \stab(S-\{s\}\cup \{r\}) = \stab(S-\{s\}\cup (W-\{s\} \cup \{r\})) = \stab(S-\{s\}) \cap \stab(W') \subseteq \stab(S-\{s\}) \cap \stab(s) = \stab(S)= \{\id\}$.  Thus $S'$ is a determining set. \medskip

Now we will show that $\stab(Y)\subset \stab(r)$ which will allow us to prove that $S'$ is minimal.\medskip

Suppose that $\beta\in \stab(Y)$.  Then since $Y$ is a block determining set for $B$, $\beta$ fixes $B$ and therefore $v$.  If $\beta(r)\ne r$ but $\beta(v)=v$ then the blockbranches at $v$ are invariant but not fixed under $\beta$.  Then there is a vertex $\tilde z$ of $T$, strictly between $\tilde x$ and $\tilde y$, that is fixed by some $\tilde \gamma$ where $\gamma$ moves $r$.  This contradicts our choice of $\alpha$ and $\tilde y$ and therefore cannot happen. Thus $\beta(r)=r$ and $\stab(Y)\subseteq \stab(r)$.  Thus $\stab(W) \subseteq \stab(Y) \subseteq \stab(r)$. \medskip

Since $S$ is minimal, $\stab(W)\subseteq \stab(r)$, and $s\in W\subset S$, we get that for all $t\in S$,  $\{\id\} \subsetneq \stab(S-\{t\}) =  \stab(S-\{s,t\}\cup W) = \stab(S-\{s,t\})\cap \stab(W) \subseteq \stab(S-\{s,t\})  \cap \stab(r) = \stab(S\cup \{r\} - \{s,t\}) = \stab(S'-\{t\})$. Thus $S'$ is also minimal.\medskip

\noindent {\bf Case 2:} Suppose that $\tilde y$ is a C-vertex. Then $H$ and $\alpha(H)$ share a cutvertex, say $v$, which tells us that there is a nontrivial isomorphism class of blockbranches at $v$.  By Condition 1, both $R$ and $S$ contain a vertex from all but one of the blockbranches of each isomorphism class of blockbranches at $v$. In particular, this is true for the, say $n$, blockbranches isomorphic to $H$ at $v$.  For $S$ denote these vertices by $s_0, s_1,\ldots, s_{n-2}$.  Since $r$ is also on one of these $n$ blockbranches, there is a subset of $n-2$ of $\{s_0,\ldots, s_{n-2}\}$ that when unioned with $r$ provides a set of vertices from all but one of the $n$ isomorphic blockbranches.  Let $s=s_0$ be a vertex that is not necessary for this and let $S'=S-\{s\}\cup\{r\}$.\medskip

Suppose that there is $\beta \in \stab(\{r, s_1,\ldots, s_{n-2} \})$ so that $\beta(s)\ne s$.  Since $\beta$ fixes each of $r, s_1,\ldots, s_{n-2}$ it leaves invariant each of their  associated blockbranches and therefore leaves invariant the remaining blockbranch in the isomorphism class of $H$ at $v$.  Thus the blockbranch, say $K$, that $s$ is on, is invariant  under $\beta$.  But since $\beta(s) \ne s$, $K$ is not fixed by $\beta$. Since $S$ is a minimal determining set, $s$ either fulfills Condition 1 of Theorem \ref{OPcriteria} for some block of $K$, or $s$ fulfills Condition 2 for some interior cutvertex of $K$.  Since $K$ and $H$ are in the same isomorphism class, $H$ must also  contain a block for which Condition 1 must be met or a cutvertex for which Condition 2 must be met. Since $R$ is a determining set, $R$ must also contain a vertex, say $t$, fulfilling the same condition.  But then $t$ also fulfills Condition 2 for $v$ and thus $r$ is unnecessary to fulfill this condition.  However, our choice of $v$ and the  minimality of $R$ prevent this situation.  Thus $\beta$ fixes $s$ and $\stab(\{r, s_1,\ldots, s_{n-2}\})\subseteq \stab(s)$.  Again using our set theoretic argument since $\{s_1,\ldots, s_{n-2}\}\subseteq S-\{s\}$  and $\stab(\{r, s_1,\ldots, s_{n-2}\})\subseteq \stab(s)$, we can show that $\stab(S')=\{\id\}$.  Similarly we can show that $\stab(\{s, s_1,\ldots, s_{n-2}\})\subseteq \stab(r)$ and therefore for all $t\in S'$, $\stab(S'-\{t\}) \ne \{\id\}$.  Thus $S'$ is a minimal determining set.\medskip

\noindent$\Longleftarrow$  \ \ Suppose that $G$ has a central block $B$ with trivial stabilizer, and that $B$ has a vertex $x$ moved by every reflection of $B$ under $\aut_G(B)$, as well as vertices $u$ and $v$ fixed by different reflections of $B$ under $\aut_G(B)$.    Then both $\{x\}$ and $\{u,v\}$ are minimal block determining sets for $B$.  However, since the pointwise stabilizer of $B$ is trivial, a block determining set for $B$ is a determining set for $G$.  Thus $\{x\}$ and $\{u,v\}$ are minimal determining sets of different size for $G$.  Thus the exchange property does not hold for determining sets in $G$.\end{proof}


\section{Open Questions}\label{sect:quest}

\begin{quest} \rm We know that determining sets do, but resolving sets do not, have the exchange property in $n$-wheels with $n\geq 8$.  Is there an infinite family of graphs in which resolving sets do, but determining sets do not, have the exchange property?\end{quest}

\begin{quest}\rm We have found precisely when determining sets have the exchange property in a given outerplanar graph.  Is there a similar characterization for resolving sets?\end{quest}

\begin{quest}\rm In which planar graphs does the exchange property hold for determining sets?  For resolving sets?  \end{quest}

\begin{quest} \rm To prove results regarding the exchange property for determining sets in trees and outerplanar graphs, we used knowledge of how automorphisms behave on these graphs.  Can we move away from detailed analyses of automorphisms and still get results on the exchange property for determining sets?  That is, is there an overarching property of a graph which tells whether the exchange property holds for determining sets?  For resolving sets? \end{quest}

\begin{quest}\rm When does the exchange property hold for distinguishing sets, locating dominating sets? \end{quest}

\bibliographystyle{plain}
\bibliography{Boutin}

\begin{thebibliography}{10}

\bibitem{AB2}
Michael~O. Albertson and Debra~L. Boutin.
\newblock Using determining sets to distinguish {K}neser graphs.
\newblock {\em Electron. J. Combin.}, 14(1):Research Paper 20 (electronic),
  2007.

\bibitem{AB1}
Michael~O. Albertson and Debra~L. Boutin.
\newblock Automorphisms and distinguishing numbers of geometric cliques.
\newblock {\em Discrete Comput. Geom.}, 39(1):778--785, 2008.

\bibitem{Ba}
L{\'a}szl{\'o} Babai.
\newblock On the complexity of canonical labeling of strongly regular graphs.
\newblock {\em SIAM J. Comput.}, 9(1):212--216, 1980.

\bibitem{B3}
Debra~L. Boutin.
\newblock The determining number of a {C}artesian product.
\newblock preprint.

\bibitem{B1}
Debra~L. Boutin.
\newblock Identifying graph automorphisms using determining sets.
\newblock {\em Electron. J. Combin.}, 13(1):Research Paper 78 (electronic),
  2006.

\bibitem{BCPZ}
Peter~S. Buczkowski, Gary Chartrand, Christopher Poisson, and Ping Zhang.
\newblock On {$k$}-dimensional graphs and their bases.
\newblock {\em Period. Math. Hungar.}, 46(1):9--15, 2003.

\bibitem{CHMPPSW}
Jos\'e C\'aceres, Carmen Hernando, Merc\`e Mora, Ignacio~M. Pelayo,
  Mar\'{\i}a~L. Puertas, Carlos Seara, and David~R. Wood.
\newblock On the metric dimension of {C}artesian products of graphs.
\newblock {\em SIAM J. Discrete Math.}, 21:423--441, 2007.

\bibitem{CEJO}
Gary Chartrand, Linda Eroh, Mark~A. Johnson, and Ortrud~R. Oellermann.
\newblock Resolvability in graphs and the metric dimension of a graph.
\newblock {\em Discrete Appl. Math.}, 105(1-3):99--113, 2000.

\bibitem{CB}
Charles~J. Colbourn and Kellogg~S. Booth.
\newblock Linear time automorphism algorithms for trees, interval graphs, and
  planar graphs.
\newblock {\em SIAM J. Comput.}, 10(1):203--225, 1981.

\bibitem{CSS}
Charles~J. Colbourn, Peter~J. Slater, and Lornal~K. Stewart.
\newblock Locating dominating sets in series parallel networks.
\newblock {\em Congr. Numer.}, 56:135--162, 1987.
\newblock Sixteenth Manitoba conference on numerical mathematics and computing
  (Winnipeg, Man., 1986).

\bibitem{DM}
John~D. Dixon and Brian Mortimer.
\newblock {\em Permutation groups}, volume 163 of {\em Graduate Texts in
  Mathematics}.
\newblock Springer-Verlag, New York, 1996.

\bibitem{EH}
David~J. Erwin and Frank Harary.
\newblock Destroying automorphisms by fixing nodes.
\newblock {\em Discrete Math.}, 306(24):3244--3252, 2006.

\bibitem{RS}
Douglas~F. Rall and Peter~J. Slater.
\newblock On location-domination numbers for certain classes of graphs.
\newblock In {\em Proceedings of the fifteenth Southeastern conference on
  combinatorics, graph theory and computing (Baton Rouge, La., 1984)},
  volume~45, pages 97--106, 1984.

\bibitem{S}
Peter~J. Slater.
\newblock Leaves of trees.
\newblock In {\em Proceedings of the Sixth Southeastern Conference on
  Combinatorics, Graph Theory, and Computing (Florida Atlantic Univ., Boca
  Raton, Fla., 1975)}, pages 549--559. Congressus Numerantium, No. XIV,
  Winnipeg, Man., 1975. Utilitas Math.

\bibitem{W}
Douglas~B. West.
\newblock {\em Introduction to graph theory}.
\newblock Prentice Hall Inc., Upper Saddle River, NJ, 1996.

\end{thebibliography}

\end{document}